\documentclass[12pt,reqno]{article}

\usepackage[usenames]{color}
\usepackage{amssymb}
\usepackage{amsmath}
\usepackage{amsthm}
\usepackage{amsfonts}
\usepackage{amscd}
\usepackage{graphicx}

\usepackage[colorlinks=true,
linkcolor=webgreen,
filecolor=webbrown,
citecolor=webgreen]{hyperref}

\definecolor{webgreen}{rgb}{0,.5,0}
\definecolor{webbrown}{rgb}{.6,0,0}

\usepackage{color}
\usepackage{fullpage}
\usepackage{float}

\usepackage{graphics}
\usepackage{latexsym}
\usepackage{epsf}
\usepackage{breakurl}

\begin{document}


\theoremstyle{plain}
\newtheorem{theorem}{Theorem}
\newtheorem{corollary}[theorem]{Corollary}
\newtheorem{lemma}[theorem]{Lemma}
\newtheorem{proposition}[theorem]{Proposition}

\theoremstyle{definition}
\newtheorem{definition}[theorem]{Definition}
\newtheorem{example}[theorem]{Example}
\newtheorem{conjecture}[theorem]{Conjecture}

\theoremstyle{remark}
\newtheorem{remark}[theorem]{Remark}

\begin{center}
\vskip 1cm{\LARGE\bf A Balanced Three-term Generalization of}
\vskip 0.3cm {\LARGE\bf Nicomachus' Identity}
\vskip 1cm
\large
Seon-Hong Kim\\
Department of Mathematics and Research Institute of Natural Science \\ 
Sookmyung Women's University\\ Seoul, 140-742\\
Korea\\
shkim17@sookmyung.ac.kr \\
\ \\
Kenneth B. Stolarsky \\ Department of Mathematics \\
University of Illinois at Urbana-Champaign\\ Urbana, IL 61801\\
USA\\
stolarsky.ken@gmail.com
\end{center}

\vskip .2 in

\begin{abstract} 
We present a generalization of the classical Nicomachus' identity for the sum of the first $n$ cubes. Unlike
previous generalizations, it has three rather than two terms, and involves not just one, but two distinct
triangular numbers, and each term is of degree $4$ in $\lfloor n/2 \rfloor$. The asymptotic behavior for large $n$ leads to continued
fractions with remarkable (but conjectural) properties. Moreover, we give a way of looking at squares of
triangular numbers that involves the square root of $11$ and show it is a limiting case of a non-obvious identity
involving truncations of the continued fraction expansion of that square root. The details involve a nonlinear
recurrence that (with appropriate initial conditions) unexpectedly produces only integers, a ``Somos-type''
phenomenon. 
    \end{abstract}
    
\section{Introduction}
We have not found formulas in the mathematical literature that involve sums of the form
\begin{equation}
\sum_{j=1}^m \bigl(a+bj+cj^2\bigr) \bigl(d+ej+fj^2\bigr),
\label{gga122}
\end{equation}
at least where all or most of the $(a, b, c, d, e, f)$ are non-zero. For $a=c=d=e=0$ and $b=f=1$, we have the classical Nicomachus' identity
\begin{equation}
\sum_{j=1}^m j^3=\Bigl(\sum_{j=1}^m j\Bigr)^2,
\label{noco1}
\end{equation}
where $T_m:=\sum_{j=1}^m j=m(m+1)/2$ is the $m$th triangular number. It is, however, possible to rewrite this in the form
$$
\sum_{j=0}^{m-1} ( 1+j )  ( 1 + 2 j + j ^ 2 )=\Bigl(\sum_{j=0}^{m-1} ( 1+j )\Bigr)^2
$$
so now only $c$ is zero. This leads immediately to a (new?) infinite sequence of Nicomachean-type identities, 
namely
$$
\sum_{j=0}^{m-1}  ( 1+a j  )  ( 1 + 2 j + a j ^ 2 ) =\Bigl(\sum_{j=0}^{m-1} ( 1+aj )\Bigr)^2.
$$
Note that this includes the well-known observation that the sum of consecutive odd numbers is a perfect square--just set the new parameter $a$ to zero. 

A different sort of generalization of (\ref{noco1}) is Liouville's identity
\begin{equation}
\sum_{d\vert n }\tau(d)^3=\Bigl( \sum_{d\vert n}\tau(d)\Bigr)^2,
\nonumber
\end{equation}
where $\tau(d)$ is the number of positive divisors of $d$. For $n$ a prime power, this becomes (\ref{noco1}). There is another generalization
\begin{equation}
    \sum_{j=1}^n \frac 1{r+2} (2j+r)\Bigl(                                                                              \begin{array}{c}
   j+r\\
   r+1 \\                                                                           \end{array}
    \Bigr)^2=\Bigl(\frac n{r+2} \Bigl(                                                                              \begin{array}{c}
   n+1+r\\
   r+1 \\                                                                           \end{array}
    \Bigr)\Bigr)^2,
    \label{genfp}
    \end{equation}
that was obtained in \cite{Ce}. We note the $1/2 \cdot 2m \cdot m^2$ pattern on the left when $r=0$, the case in which this reduces to (\ref{noco1}). An elegant $q$-identity generalization of (\ref{noco1}) is Warnaar's $q$-analogue (see \cite{Wa}), namely
\begin{equation}
\sum_{j=1}^n q^{2n-2j}\frac {1-q^{2j}}{1-q^2}\Bigl(\frac {1-q^j}{1-q}\Bigr)^2 =\frac 1{(1+q)^2} \Bigl(\frac {1-q^n}{1-q} \frac {1-q^{n+1}}{1-q}\Bigr)^2.
\nonumber
\end{equation}
Compared to (\ref{genfp}), this has $1\cdot m \cdot m^2$ pattern on the left as $q \rightarrow 1$, the case in which this reduces to (\ref{noco1}). We now recall that the fact that $n^2$ is the sum of consecutive odd numbers can be shown by taking a square array of $n^2$ unit squares and expressing it as a union of $n$ gnomons having $1$, $3$, $5$, $\ldots$, $2 n - 1$ unit squares respectively. More than $1000$ years ago Arab mathematicians had an analogous proof for (\ref{noco1}). See \cite [p. 32 and note 23 on p. 187] {Ot}. Such geometric proofs can be viewed as applications of ``figurate numbers'' and are sometimes known as ``proofs without words''. These have accumulated sporadically over the centuries. To organize them and produce them in a systematic way seems to us to be the essence of \cite{Ma}. Identity (\ref{noco1})
appears there as Example 11, p. 15. It begins by writing  the integer summand $i ^ 3$ as $i \cdot i ^ 2$ and then handles each of these two factors in a different way. We now note that in our expression (\ref{gga122}), there is no $j ^ 3$, but when $c = 0$, the lead terms of the summand factors are $j$ and $j ^ 2$. Altogether we have thus found the $j \cdot j ^ 2$ pattern in four distinct cases.

In Section~\ref{mthm} and Section~\ref{sect2}, we state and prove a new generalization of (\ref{noco1}) involving a parameter $x$ in which the expression on the right that is squared involves two triangular numbers. The results are (\ref{rt1}) and (\ref{rt2}). These incorporate a sum of the form (\ref{gga122}), with $c=0$, into a three term identity that generalizes (\ref{noco1}). 

We have not previously seen product expressions such as $x P-x^2+x^3$ in (\ref{rt1}) and (\ref{rt2}) in common identities, but they do have a role to play. Also, a bit of the fascination of the Nicomachean-type identity is from the fact that the cube sum is given by a polynomial with the (statistically rare) property of having two zeros of multiplicity two. We ask: when do linear combinations of power sum formulas have zeros of multiplicity at least two? The expression (\ref{gga122}) is obviously a linear combination of power sums. The answer to the question here is the set of points in the 6-dimensional $(a,b,c,d,e,f)$ space, where the discriminant $D$ with respect to $m$ of the expression is $0$. Computer algebra shows, alas, that $D$ has $1114$ terms with large coefficients. Such unwieldy expressions are a common obstruction to progress in this area as many of the unresolved conjectures in \cite{Ki} indicate. However, in this case computer algebra does show that $D$ has a factor $F$ given by
$$F=5a(6d+3e+f)+5b(3d+e)+c(5d-f).
$$
The Nicomachean case corresponds to $a=b=d=f=0$ which satisfies $F=0$. If
\begin{equation}
(a, b, c)=\Bigl(1, \,\frac {- 1 + \sqrt {11}}2, \,\frac {- 1 - \sqrt {11}}2\Bigr)\,\,\, \text{and} \,\,\, (d,e,f)=(-a,-c,-b),
\label{abccb}
\end{equation}
then (\ref{gga122}) becomes
\begin{equation}
( 2 m + 3 ) T_{ m - 1 }^2   ,
\label{gdmoide}
\end{equation}
almost the classical Nicomachus' identity. We suspect it is new, and a good illustration of Littlewood's dictum, that identities are trivial, but only after they have been discovered. 
Moreover, we show it is a limiting case of a more general identity in Section~\ref{moge11}. 

The coefficients in (\ref{abccb}) do satisfy $F = 0$. For another example, $(a,d,e,f)=(0,1,-3,5)$ gives $F=0$ and
\begin{align*}
&\sum_{j=1}^m (bj+cj^2)(1-3j+5j^2)\\
=& \frac 14 m^2 (1+m)\bigl(b-c+(5b+3c)m+4cm^2 \bigr).
\end{align*}
We observe that in the classical Nicomachus' identity, the multiple zeros are a distance $1$ apart. If a sum of the form (\ref{gga122}) has two multiple zeros, how far apart could they be? We leave this question open, but the following identities with $c=0$ shed some light on this:
\begin{equation}
\sum_{j=1}^m (1-bj)(1-2j+bj^2)=-\frac 14 m^2 (b-2+bm)^2,
\nonumber
\end{equation}
\begin{equation}
\sum_{j=1}^m (1-bj)(1+(b-2)j)\bigl(1+b(b-2)j \bigr)=-\frac 14 m^2 \bigl( (b-1)^2+1+m((b-1)^2-1)\bigr)^2.
\nonumber
\end{equation}

Considerations of natural generalization of (\ref{gga122}) would lead into a combination of Partition theory (see \cite{An} for an introduction) and the study of power sums (see e.g. \cite{Kn}). This is beyond the scope of this paper. Recently, by introducing translation parameters, adjoining variable parameters, and forming sequential products, we created many general Nicomachean-type identities in \cite{Ki2}. The paper \cite{Bh} provides sophisticated elliptic extensions of several identities including (\ref{noco1}). For example, a special case of an identity in \cite{Bh} yields Cigler's $q$-analogue of (\ref{noco1}) (See \cite{Ci}), namely
$$
\sum_{j=1}^m q^{\binom{m+1}{2}-\binom{j+1}{2}}\frac {1-q^j}{1-q}\frac {1-q^{j^2}}{1-q}=\Bigl(\frac {1-q^{\binom{m+1}2}}{1-q}\Bigr)^2.
$$

There is more to say about the interaction of discriminants and identities of Nicomachean type. See \cite{Ki}.

\section{Theorem~\ref{thj1}}\label{mthm}

Balance in significant three-term identities is not pre-ordained. For example, in the Pythagorean triplet identity
$$
(2m-1)^2+(2m(m-1))^2=(m^2+(m-1)^2)^2,
$$
the first term grows like $4m^2$ while the others grow like $4m^4$.

The first part of Theorem~\ref{thj1} below generalizes Nicomachus' identity for the odd termination sequence $\{1, 2, \cdots, 2m+1\}$. The second part deals with the case of even termination. The present identities in (\ref{rt1}) and (\ref{rt2}) below are balanced in $m$. The $L$, $R$ and $xP-x^2+x^3$ terms all have degree $4$ in $m$. Moreover, the coefficients of $m^4$ in these three cases factor nicely. They are, respectively,
\begin{equation}
\frac 14 (x+4)\bigl((x+1)^2+3\bigr), \,\,\, \frac 14(x+4)^2, \,\,\, \text{and}\,\,\, \frac 14 x(x+1)(x+4).
\label{betx1}
\end{equation}
This implies that for $m$ very large, each such expression has a root very near to $-4$. We explore this further in Remark~\ref{cfrac5}.

We observe that the respective coefficients of $x^3$ are
$$
\Big\{ 1+\Bigl(\frac {m(m+1)}2\Bigr)^2, \,0,  \,1+\Bigl(\frac {m(m+1)}2\Bigr)^2\Big\}.
$$
So the squares of the triangular numbers are hidden in the first and third terms, rather than being explicit as in the middle term $R$.

\begin{theorem} Let $m\geq 2$ be an integer. Let
\begin{align*}
L_{m,o}(x)&=\sum_{j=1}^m (2 j -1 ) ^ 3   + x ^ 3 + \sum_{j=2}^m ( ( 2 j - 2 ) + ( j - 1 ) x ) ^ 3 ,\\
R_{m,o}(x) &= (  T_{2 m - 1} + ( 1 + T_{ m - 1}) x )^ 2.
\end{align*}
Then we have
\begin{equation}
L_{m+1,o}(x)-R_{m+1,o}(x)=xP_{m,o}(x)-x^2+x^3,
\label{rt1}
\end{equation}
where
$$
P_{m,o}(x)=\sum_{j=0}^m \Bigl( ( 4 j + 1 + j x ) \bigl(  ( j + 1 ) ( j - 2 ) + j ^ 2  x \bigr) \Bigr).
$$
Let
\begin{align*}
L_{m,e}(x)&=\sum_{j=1}^m (2 j -1 ) ^ 3   + x ^ 3 + \sum_{j=1}^m (  2 j   + j x ) ^ 3 ,\\
R_{m,e}(x) &= (  T_{2 m } + ( 1 + T_{ m }) x )^ 2.
\end{align*}
Then
\begin{equation}
L_{m,e}(x)-R_{m,e}(x)=xP_{m,e}(x)-x^2+x^3,
\label{rt2}
\end{equation}
where $P_{m,e}(x)$ is the polynomial obtained by substituting $-m-1$ for $m$ in the polynomial expansion of $P_{m,o}(x)$.
\label{thj1}
\end{theorem}

\begin{remark} When $x=0$, equations (\ref{rt1}) and (\ref{rt2}) become
$$
\sum_{j=1}^{2m+1}j^3=T_{2m+1}^2\,\,\, \text{and}\,\,\, \sum_{j=1}^{2m}j^3=T_{2m}^2,
$$
respectively, which together are Nicomachus' identity. For $m\geq 3$, the quadratic polynomials $P_{m,o}(x)$ and $P_{m,e}(x)$ have positive integer coefficients. For $x$ very large, the summand of $P_{m,o}(x)$ behaves like $j ^ 3$, and this arises as the product of $j$ with $j ^ 2$, the pattern we have noted among two of the previously published generalizations. The leading coefficient of $P_{m,o}(x)$ is $T_m^2$, and so both the left and right sides in (\ref{rt1}) behave like $T_m^2+1$ as $x\rightarrow \infty$. Since $T_{-m-1}^2=T_m^2$, the same behavior as $x\rightarrow \infty$ as mentioned in (\ref{rt1}) holds for (\ref{rt2}).
\end{remark}

\begin{remark}  Equations (\ref{rt1}) and (\ref{rt2}) involves the expression $- x ^ 2+x^3$. This is a very natural expression in the study of Nicomachean-type identities. Take the sequence $\{ 1 , 2 , 3 , \ldots , n \}$ and replace one of the 
integers by a variable $x$. Then calculate and simplify the square of the sum minus the
sum of the cubes, The constant coefficient and the coefficient of $x$ will depend upon the
value of $n$, but the remaining part will always be $- x ^ 2+x ^ 3$.
\end{remark}

\begin{remark} For large $m$, each of the three terms $L$, $R$ and $xP-x^2+x^3$ in (\ref{rt1}) and (\ref{rt2}) has a root $r$ near $-4$. The continued fractions of the negatives of these roots are of interest. For the odd termination case, one expects the second partial quotients to be large, but one soon finds yet another large partial quotient. We shall examine them for a sequence of $m$-values tending to infinity that are chosen for convenience, namely $m=10^{2(2n+1)}$, $n=1,2,3,\ldots$. For large $n$, we conjecture the following: for the $L_{m+1,o}$ term, the $R_{m+1,o}$ term, the $xP_{m,o}-x^2+x^3$ term, the negatives of the roots have continued fraction representations
\begin{align*}
&\Bigl[4\, ; \, 75\cdot 10^{4n}, \, 3, \, 4, \, 1,\, \frac {(15\cdot 10^{2n})^2-23}{101}, \, \ldots \Bigr],\\
&\Bigl[4\, ; \, 50\cdot 10^{4n}+2, \, \left\lfloor\frac { 10^{4n+2}}{7}\right\rfloor-1, \, \ldots\Bigr],\\
&\Bigl[4\, ; \, 75\cdot 10^{4n}, \, 3, \, \left\lceil\frac {6\cdot 10^{4n+2}+100}{35}\right\rceil, \, \ldots\Bigr],
\end{align*}
respectively. Here the formula for the $R_{m+1,o}$ term gives the largest value. In the formulas for $L_{m+1,o}$ and $xP_{m,o}-x^2+x^3$, the value for $L_{m+1,o}$ is larger. We also observe that in each case, the ratio of the second partial quotient to the next large partial quotient is very close to $101/3$, $7/2$ and $35/8$ when $n$ is large, respectively.

There are similar phenomena for the even termination case. 
\label{cfrac5}
\end{remark}

\section{Proof of Theorem~\ref{thj1}}\label{sect2}
The results in Theorem~\ref{thj1} can be given a purely mechanical proof using the (many centuries old) formulas for the sums of the $n$-th powers of consecutive integers, $1 \leq n \leq 3$. However, there is a nice way of presenting
such a proof when it is ``half-way'' done. Compute $L_{m+1,o}$, $R_{m+1,o}$, $x P_{m,o}$ as polynomials in $m$ and $x$. Then create the $4 \times 3$ matrix $M$ whose columns are the coefficients of those three viewed as polynomials in $x$. Then $M$ is
\begin{equation}
\begin{bmatrix}
     T_{2m+1}^2       & (1 + m)^2 (1 + 2 m)^2  & 0\\
   12  T_ m^2       & (1 + m) (1 + 2 m) (2 + m + m^2) & (1 + m) (2 + m) (-1 - 2 m + m^2) \\
   6  T_ m^2      &  \frac 14 (2 + m + m^2)^2 & \frac 14 m (1 + m) (-4 + 5 m + 5 m^2)\\
         1 +  T_ m^2  &   0 &  T_ m^2
\end{bmatrix}.
\label{mat1}
\end{equation}
Now the product of the matrix $M$ and the transpose of $[ 1 , - 1 , - 1 ]$ is simply the transpose of $[ 0 , 0 , -1 , 1]$, and this corresponds to $- x ^ 2 + x ^ 3$, thus completing
the proof of the first part of the theorem. One can find a corresponding matrix for the second part of the theorem. In fact, for
$$
J=-(2m+1)^3, \quad K=(2m+1)(2+m+m^2),
$$
if the matrix
\begin{equation}
\begin{bmatrix}
    J      &J  & 0\\
          0& -K & K\\
        0&  0 & 0\\
        0  &   0 &  0
\end{bmatrix}
\label{mat3}
\end{equation}
multiplies the transpose of $[1 , -1 , -1  ]$, we clearly get the zero vector. Now calculations show that if we add the matrix (\ref{mat3}) to the matrix $M$ that arises in the proof of the first part, we get the matrix that arises in the proof of the second part.
\par
By examining both matrices, we are led to a ``double $a = b + c$
equation'' of interest, namely
\begin{align*}
12T_ m^2&=(1 + m) (1 + 2 m) (2 + m + m^2) + (1 + m) (2 + m) (-1 - 2 m + m^2)\\
&= m (1 + 2 m) (2 + m + m^2) + (-1 + m) m (2 +4 m + m^2),
\end{align*}
where $a$, $b$ and $c$ are all very composite. This may be of some interest in the study of polynomial analogues of ``abc'' type problems. We also have
\begin{equation}
7m^3 (1+m)^3=(1+m)(1+2m)(1+3m)A_3(m)+(1+m)(2+m)(3+m)B_3(m),
\nonumber
\end{equation}
where $A_3(m)$ and $B_3(m)$ are monic polynomials of degree $3$. Moreover there are infinitely many pairs of $(A_3(m), B_3(m))$ of such polynomials, with rational coefficients.

\begin{remark} In (\ref{mat1}), we note that if the $j$-th row $( 1 \leq j \leq 4 )$ of the $4 \times 3$ matrix is removed, we have
in each case a $3\times 3$ matrix. The determinants of these four matrices are notably composite.
They are
\begin{align*}
&\frac 14 (1 + m) (2 + m + m^2) (4 + 8 m + m^2 + 19 m^3 + 29 m^4 + 11 m^5),\\
&\frac 12 m (1 + m)^3 (1 + 2 m)^2 (-2 + 3 m + 3 m^2),\\
&(1 + m)^3 (2 + m) (1 + 2 m)^2 (-1 - 2 m + m^2),\\
&-(1 + m)^3 (2 + m) (1 + 2 m)^2 (-1 - 2 m + m^2).
\end{align*}
Note that the last two are negatives of each other. One could also interchange the roles of $x$ and $m$ and finish the proof by multiplying a
$3 \times 5$ matrix on the left by $[ 1 , -1 , -1 ]$. This proof also provides the representation
$$
- x ^ 2 + x ^ 3 =
           (1 + x ) (1 - x + x ^ 2 ) -  ( 1 + x ) ^ 2 - ( - 2 x ) .
$$
Also, if the $3$rd and $4$th columns are removed, the resulting $3 \times 3$ matrix has the rather composite determinant
$$
\frac 14 x^3(-1+ x  )   (4+ x )  ( 34 + 24 x + 9 x ^ 2 + x ^ 3 ) .
$$
\end{remark}

\section{The ``$\sqrt {11}$ identity'' as a limiting case}\label{moge11}

A type (\ref{gga122}) formula, namely
\begin{equation}
\frac 1{2m+3} \sum_{j=1}^m \bigl(1+bj+cj^2\bigr)\bigl(-1-cj-bj^2\bigr)=T_{m-1}^2
\label{2m35}
\end{equation}
with $b=\bigl(-1+\sqrt{11}\bigr)/2$ and $c=\bigl(-1-\sqrt{11}\bigr)/2$, gives a new way of looking at squares of triangular numbers. This is almost the classical Nicomachus' identity. 

The continued fraction of $(-1+\sqrt {11})/2$ is
$$
[1; \overline{6, 3}],
$$
and a convergent of this is
\begin{equation}
\alpha_k:=[1; \underbrace{6,3,6,3,\ldots, 6,3}_{2k-2 \,\,\text{numbers}}].
\label{alk1012}
\end{equation}
Then
\begin{equation}
\alpha_{k+1}=1+\frac 1{6+\frac 1{3+(\alpha_k-1)}}=\frac {15+7\alpha_k}{13+6\alpha_k}.
\label{alk00}
\end{equation}

The positive integers $u_k$ below in (\ref{110hj}) will be of special interest in connection with $\alpha_k$. Set $A = ( 2 \cdot 3 \cdot 5 )^2 = 900$. We note that the discriminant of $1 - 398 x + x^2$ is $A \cdot 2^4 \cdot 11$ and every term in the power series expansion of Proposition~\ref{prp3} is congruent to $1$ modulo $A$. Also, the discriminant of $1 + 502 x + x ^ 2$ has only small prime factors, namely $A\cdot 2 ^ 3 \cdot 5 \cdot 7$.

\begin{proposition} Let $u_k$ be the coefficient of $x^k$ in the power series expansion of
\begin{equation}
\frac {1 + 502 x + x ^ 2 }{( 1 - x ) ( 1 - 398 x + x ^ 2 )}=1 + 901 x + 359101 x^2 + 142921801 x^3 + \cdots.
\nonumber
\end{equation}
Then
\begin{equation}
u_{k}=\frac{1}{44} \bigl(5 \bigl(10-3 \sqrt{11}\bigr)^{2
   k-1}+5\bigl(10+3 \sqrt{11}\bigr)^{2 k-1}-56\bigr)
\label{110hj}
\end{equation}
and
\begin{equation}
u_{k+1}=252 + 199 u_k +
 30 \sqrt {(3 + 2 u_k) (23 + 22 u_k)}.
\label{uk00}
\end{equation}
\label{prp3}
\end{proposition}
\begin{proof} Using (\ref{110hj}) and computer algebra, we find that
\begin{align*}
&\sum_{k=1}^n u_k x^{k-1}\\
 =& \frac 1{44 (x^3-399
   x^2+399 x-1)}\\
 &\Bigl(44 u_n
   x^{n+2} -4 \bigl(25 \bigl(10-3 \sqrt{11}\bigr)^{2 n}+25
   \bigl(10+3 \sqrt{11}\bigr)^{2 n}-5572\bigr)
   x^{n+1}+44u_{n+1} x^n\\
& -44 x^2-22088 x-44\Bigr).
\end{align*}
For $x\in (-x_0,x_0)$ for sufficiently small $x_0>0$, the first three terms of the second factor in the above tends to $0$ as $n\rightarrow \infty$. So if $x\in (-x_0,x_0)$,
$$
\lim_{n\rightarrow \infty}\sum_{k=1}^n u_k x^{k-1}=
\frac {-44 x^2-22088 x-44}{44 (x^3-399
   x^2+399 x-1)}=\frac {1 + 502 x + x ^ 2 }{( 1 - x ) ( 1 - 398 x + x ^ 2 )}.
   $$
The recurrence relation in (\ref{uk00}) follows from (\ref{110hj}).
\end{proof}

The number $\alpha_k$ in (\ref{alk1012}) can be represented by $u_k$.

\begin{proposition} For $k\geq 1$
\begin{equation}
\alpha_k=\frac 12 \Bigl( -1+\sqrt{11-\frac {10}{3 + 2 u_k}}\Bigr).
\label{iihh07}
\end{equation}
\end{proposition}
\begin{proof} Let  
\begin{equation}
\beta_k=\frac 12 \Bigl( -1+\sqrt{11-\frac {10}{3 + 2 u_k}}\Bigr).
\nonumber
\end{equation}
Then by computer algebra, (\ref{uk00}) implies that
$$
\beta_{k+1}-\frac {15+7\beta_k}{13+6\beta_k}=0.
$$
Since $\beta_1=1$, we have $\beta_k=\alpha_k$ for all $k\geq 1$.
\end{proof}

\begin{theorem} For
\begin{equation}
(a, b, c)=\bigl(1, \,\alpha_k, \,-(1+\alpha_k)\bigr)\,\,\, \text{and} \,\,\, (d,e,f)=(-a,-c,-b),
\label{ggtt58}
\end{equation}
we have
$$
\sum_{j=1}^m (a+bj+cj^2)(d+ej+fj^2)=T_{ m - 1 }\Bigl(  \frac {(1+ m ) (2- 3 m ^ 2) }{
  3  ( 3 + 2 u_k ) }  + ( 3+2 m ) T_{ m - 1 }\Bigr).
$$
\end{theorem}

\begin{proof} Straightforward computations using (\ref{iihh07}) and (\ref{ggtt58}) yield
\begin{align*}
&(a+bj+cj^2)(d+ej+fj^2)\\
=&a d + (b d + a e) j + (c d + b e + a f) j^2 + (c e + b f) j^3 + 
 c f j^4\\
=& -1 + j + (1 + \alpha_k +  \alpha_k^2) j^2 - (1 + 2\alpha_k+2\alpha_k^2) j^3 + \alpha_k(1 + 
   \alpha_k) j^4\\
=&-1 + j +\frac {8+7u_k}{3+2u_k} j^2-\frac {13+12u_k}{3+2u_k}j^3+\frac {5(1+u_k)}{3+2u_k} j^4.   
\end{align*}
Using the above, power sum formulas and simplifying them, we obtain
\begin{align*}
&\sum_{j=1}^m (a+bj+cj^2)(d+ej+fj^2)\\
 =&\frac {(-1 + m) m}2\cdot \frac {4 - 23 m + 3 m^2 + 12 m^3 + 6 m (-3 + m + 2 m^2) u_k}{6(3+2u_k)}\\
=& T_{ m - 1 }\cdot \frac { 2(1 + m) (2 - 3 m^2)+3(-1 + m) m (3 + 2 m)(3+2u_k)}{6(3+2u_k)}\\
 =&T_{ m - 1 }\Bigl(  \frac {(1+ m ) (2- 3 m ^ 2) }{
  3  ( 3 + 2 u_k ) }  + ( 3+2 m ) T_{ m - 1 }\Bigr).
\end{align*}
\end{proof}

\begin{remark} The identity (\ref{2m35}) is simply the $k \rightarrow \infty$ limiting case of the above theorem
since $u_k \rightarrow \infty$ as $k \rightarrow \infty$.
\end{remark}

\begin{remark} Congruences such as the above for modulus $900$, mentioned just before the statement of Proposition~\ref{prp3}, can be ``artificially'' produced.
Say a ``congruent to $1$ mod $m$'' congruence is desired with m ``somewhat large''. Set
$$
f( m , x ) = \frac {m x}{(  1 - x ) ( 1 - 123 x + x ^ 2 )},
$$
where $123$ has been chosen simply because the discriminant of $1 - 123 x + x ^ 2$, namely
$5 ^ 3 \cdot 11^2$, has only ``small'' prime factors. Next, note that the discriminant of $1 + 130 x + x ^ 2$, namely $2 ^ 9 \cdot 3 \cdot 11$, also has only ``small'' prime factors. Now
$$
f ( m , x )  + \frac 1{ 1 - x } =  \frac { 1 + ( m - 123 ) x + x ^ 2}{ ( 1 - x ) ( 1 - 123 x + x ^ 2 )}
$$
and $m - 123 = 130$ when $m = 253$. Thus all power series coefficients of the above rational
expression for this $m$, which now has $1 + 130 x + x ^ 2$ in the numerator, are congruent to $1$ mod $253$. However, in the present work the congruence came unbidden with no such construction employed.
\end{remark}

\bigskip
\hrule
\bigskip

\noindent 2020 {\it Mathematics Subject Classification}:
Primary 11B83; Secondary 11D25.

\noindent \emph{Keywords: } sum of cubes, Nicomachus' identity, Nicomachean-type identity, triangular number, continued fraction, Somos-type phenomenon, total factorization number

\bigskip
\hrule
\bigskip

\enddocument